\documentclass[notitlepage,leqno,11pt]{amsart}
\usepackage{amsmath,amssymb,amsbsy,amsfonts,amsthm,latexsym,
            amsopn,amstext,amsxtra,euscript,amscd}
\usepackage{hyperref}

\begin{document}
\bibliographystyle{plain}
\newfont{\teneufm}{eufm10}
\newfont{\seveneufm}{eufm7}
\newfont{\fiveeufm}{eufm5}
%
%
\newfam\eufmfam
              \textfont\eufmfam=\teneufm \scriptfont\eufmfam=\seveneufm
              \scriptscriptfont\eufmfam=\fiveeufm
\def\bbbr{{\rm I\!R}}
\def\bbbm{{\rm I\!M}}
\def\bbbn{{\rm I\!N}}
\def\bbbf{{\rm I\!F}}
\def\bbbh{{\rm I\!H}}
\def\bbbk{{\rm I\!K}}
\def\bbbp{{\rm I\!P}}
\def\bbbone{{\mathchoice {\rm 1\mskip-4mu l} {\rm 1\mskip-4mu l}
{\rm 1\mskip-4.5mu l} {\rm 1\mskip-5mu l}}}
\def\bbbc{{\mathchoice {\setbox0=\hbox{$\displaystyle\rm C$}\hbox{\hbox
to0pt{\kern0.4\wd0\vrule height0.9\ht0\hss}\box0}}
{\setbox0=\hbox{$\textstyle\rm C$}\hbox{\hbox
to0pt{\kern0.4\wd0\vrule height0.9\ht0\hss}\box0}}
{\setbox0=\hbox{$\scriptstyle\rm C$}\hbox{\hbox
to0pt{\kern0.4\wd0\vrule height0.9\ht0\hss}\box0}}
{\setbox0=\hbox{$\scriptscriptstyle\rm C$}\hbox{\hbox
to0pt{\kern0.4\wd0\vrule height0.9\ht0\hss}\box0}}}}
\def\bbbq{{\mathchoice {\setbox0=\hbox{$\displaystyle\rm
Q$}\hbox{\raise
0.15\ht0\hbox to0pt{\kern0.4\wd0\vrule height0.8\ht0\hss}\box0}}
{\setbox0=\hbox{$\textstyle\rm Q$}\hbox{\raise
0.15\ht0\hbox to0pt{\kern0.4\wd0\vrule height0.8\ht0\hss}\box0}}
{\setbox0=\hbox{$\scriptstyle\rm Q$}\hbox{\raise
0.15\ht0\hbox to0pt{\kern0.4\wd0\vrule height0.7\ht0\hss}\box0}}
{\setbox0=\hbox{$\scriptscriptstyle\rm Q$}\hbox{\raise
0.15\ht0\hbox to0pt{\kern0.4\wd0\vrule height0.7\ht0\hss}\box0}}}}
\def\bbbt{{\mathchoice {\setbox0=\hbox{$\displaystyle\rm.
T$}\hbox{\hbox to0pt{\kern0.3\wd0\vrule height0.9\ht0\hss}\box0}}
{\setbox0=\hbox{$\textstyle\rm T$}\hbox{\hbox
to0pt{\kern0.3\wd0\vrule height0.9\ht0\hss}\box0}}
{\setbox0=\hbox{$\scriptstyle\rm T$}\hbox{\hbox
to0pt{\kern0.3\wd0\vrule height0.9\ht0\hss}\box0}}
{\setbox0=\hbox{$\scriptscriptstyle\rm T$}\hbox{\hbox
to0pt{\kern0.3\wd0\vrule height0.9\ht0\hss}\box0}}}}
\def\bbbs{{\mathchoice
{\setbox0=\hbox{$\displaystyle     \rm S$}\hbox{\raise0.5\ht0\hbox
to0pt{\kern0.35\wd0\vrule height0.45\ht0\hss}\hbox
to0pt{\kern0.55\wd0\vrule height0.5\ht0\hss}\box0}}
{\setbox0=\hbox{$\textstyle        \rm S$}\hbox{\raise0.5\ht0\hbox
to0pt{\kern0.35\wd0\vrule height0.45\ht0\hss}\hbox
to0pt{\kern0.55\wd0\vrule height0.5\ht0\hss}\box0}}
{\setbox0=\hbox{$\scriptstyle      \rm S$}\hbox{\raise0.5\ht0\hbox
to0pt{\kern0.35\wd0\vrule height0.45\ht0\hss}\raise0.05\ht0\hbox
to0pt{\kern0.5\wd0\vrule height0.45\ht0\hss}\box0}}
{\setbox0=\hbox{$\scriptscriptstyle\rm S$}\hbox{\raise0.5\ht0\hbox
to0pt{\kern0.4\wd0\vrule height0.45\ht0\hss}\raise0.05\ht0\hbox
to0pt{\kern0.55\wd0\vrule height0.45\ht0\hss}\box0}}}}
\def\bbbz{{\mathchoice {\hbox{$\sf\textstyle Z\kern-0.4em Z$}}
{\hbox{$\sf\textstyle Z\kern-0.4em Z$}}
{\hbox{$\sf\scriptstyle Z\kern-0.3em Z$}}
{\hbox{$\sf\scriptscriptstyle Z\kern-0.2em Z$}}}}
\def\ts{\thinspace}

\newtheorem{theorem}{Theorem}
\newtheorem{lemma}[theorem]{Lemma}
\newtheorem{claim}[theorem]{Claim}
\newtheorem{cor}[theorem]{Corollary}
\newtheorem{prop}[theorem]{Proposition}
\newtheorem{definition}[theorem]{Definition}
\newtheorem{remark}[theorem]{Remark}
\newtheorem{question}[theorem]{Open Question}
\newtheorem{example}[theorem]{Example}

\def\qed{\ifmmode
\squareforqed\else{\unskip\nobreak\hfil
\penalty50\hskip1em\null\nobreak\hfil\squareforqed
\parfillskip=0pt\finalhyphendemerits=0\endgraf}\fi}

\def\squareforqed{\hbox{\rlap{$\sqcap$}$\sqcup$}}

\def \C {{\mathbb C}}
\def \F {{\mathbb F}}
\def \L {{\mathbb L}}
\def \K {{\mathbb K}}
\def \Q {{\mathbb Q}}
\def \Z {{\mathbb Z}}
\def\cA{{\mathcal A}}
\def\cB{{\mathcal B}}
\def\cC{{\mathcal C}}
\def\cD{{\mathcal D}}
\def\cE{{\mathcal E}}
\def\cF{{\mathcal F}}
\def\cG{{\mathcal G}}
\def\cH{{\mathcal H}}
\def\cI{{\mathcal I}}
\def\cJ{{\mathcal J}}
\def\cK{{\mathcal K}}
\def\cL{{\mathcal L}}
\def\cM{{\mathcal M}}
\def\cN{{\mathcal N}}
\def\cO{{\mathcal O}}
\def\cP{{\mathcal P}}
\def\cQ{{\mathcal Q}}
\def\cR{{\mathcal R}}
\def\cS{{\mathcal S}}
\def\cT{{\mathcal T}}
\def\cU{{\mathcal U}}
\def\cV{{\mathcal V}}
\def\cW{{\mathcal W}}
\def\cX{{\mathcal X}}
\def\cY{{\mathcal Y}}
\def\cZ{{\mathcal Z}}
\newcommand{\rmod}[1]{\: \mbox{mod}\: #1}

\def\tcN{\cN^\mathbf{c}}
\def\F{\mathbb F}
\def\Tr{\operatorname{Tr}}
\def\mand{\qquad \mbox{and} \qquad}
\renewcommand{\vec}[1]{\mathbf{#1}}
\def\eqref#1{(\ref{#1})}
\newcommand{\ignore}[1]{}
\hyphenation{re-pub-lished}
\parskip 1.5 mm
\def\lln{{\mathrm Lnln}}
\def\Res{\mathrm{Res}\,}
\def\F{{\bbbf}}
\def\Fp{\F_p}
\def\fp{\Fp^*}
\def\Fq{\F_q}
\def\ff{\F_2}
\def\ffn{\F_{2^n}}
\def\K{{\bbbk}}
\def \Z{{\bbbz}}
\def \N{{\bbbn}}
\def\Q{{\bbbq}}
\def \R{{\bbbr}}
\def \P{{\bbbp}}
\def\Zm{\Z_m}
\def \Um{{\mathcal U}_m}
\def \Bf{\frak B}
\def\Km{\cK_\mu}
\def\va {{\mathbf a}}
\def \vb {{\mathbf b}}
\def \vc {{\mathbf c}}
\def\vx{{\mathbf x}}
\def \vr {{\mathbf r}}
\def \vv {{\mathbf v}}
\def\vu{{\mathbf u}}
\def \vw{{\mathbf w}}
\def \vz {{\mathbfz}}
\def\\{\cr}
\def\({\left(}
\def\){\right)}
\def\fl#1{\left\lfloor#1\right\rfloor}
\def\rf#1{\left\lceil#1\right\rceil}
\def\flq#1{{\left\lfloor#1\right\rfloor}_q}
\def\flp#1{{\left\lfloor#1\right\rfloor}_p}
\def\flm#1{{\left\lfloor#1\right\rfloor}_m}
\def\Al{{\sl Alice}}
\def\Bob{{\sl Bob}}
\def\Or{{\mathcal O}}
\def\inv#1{\mbox{\rm{inv}}\,#1}
\def\invM#1{\mbox{\rm{inv}}_M\,#1}
\def\invp#1{\mbox{\rm{inv}}_p\,#1}
\def\Ln#1{\mbox{\rm{Ln}}\,#1}
\def \nd {\,|\hspace{-1.2mm}/\,}
\def\ord{\mu}
\def\E{\mathbf{E}}
\def\Cl{{\mathrm {Cl}}}
\def\epp{\mbox{\bf{e}}_{p-1}}
\def\ep{\mbox{\bf{e}}_p}
\def\eq{\mbox{\bf{e}}_q}
\def\bm{\bf{m}}
\newcommand{\floor}[1]{\lfloor {#1} \rfloor}
\newcommand{\comm}[1]{\marginpar{
\vskip-\baselineskip
\raggedright\footnotesize
\itshape\hrule\smallskip#1\par\smallskip\hrule}}
\def\rem{{\mathrm{\,rem\,}}}
\def\dist {{\mathrm{\,dist\,}}}
\def\etal{{\it et al.}}
\def\ie{{\it i.e. }}
\def\veps{{\varepsilon}}
\def\eps{{\eta}}
\def\ind#1{{\mathrm {ind}}\,#1}
               \def \MSB{{\mathrm{MSB}}}
\newcommand{\abs}[1]{\left| #1 \right|}

\title{ Pell Numbers with Lehmer property }
\subjclass[2010]{Primary 11B39; Secondary A25}
\keywords{Pell numbers, Euler function}
\author{
{\sc Bernadette~Faye\quad Florian~Luca}
}
\address{
Ecole Doctorale de Mathematiques et d'Informatique \newline
Universit\'e Cheikh Anta Diop de Dakar \newline
BP 5005, Dakar Fann, Senegal \newline
\and \newline
School of Mathematics, University of the Witwatersrand \newline
Private Bag X3, Wits 2050, South Africa 
}
\email{bernadette@aims-senegal.org} 
\address{
School of Mathematics, University of the Witwatersrand \newline
Private Bag X3, Wits 2050, South Africa \newline
}
\email{Florian.Luca@wits.ac.za}

\maketitle

\begin{abstract} In this paper, we prove that there is no number with the Lehmer property in the sequence of Pell numbers.
\end{abstract}

\section{Introduction}
Let $\phi(n)$ be the Euler function of a positive  integer $n$. A composite integer has the Lehmer property if $\phi(n)\mid n-1$. Lehmer \cite{Leh} conjectured that there are no such numbers. For a positive integer $m$, we write $\omega(m)$ for the number of  distinct prime divisors of $m$. Lehmer proved that if $N$ has his property, then $\omega(N)\geq 7$. This result has been improved by Cohen and Hagis  in \cite{coh} to $\omega(N)\geq 14.$ The current record $\omega(N)\geq 15$ is due to Renze \cite{joh}. In case where $3\mid N$, Bursci et al. \cite{PB} proved that $\omega(N)\geq 40\cdot 10^{6}$ and $N>10^{36\cdot 10^{7}}.$

Many results concerning this problem can be found in the litterature (see \cite{FL2}, \cite{CP}). Not succeeding in proving that there are no numbers with the Lehmer property, some researchers concentrated on proving that there are no numbers with the Lehmer property in certain interesting subsequences of positive integers like in the Fibonacci sequence $\{F_n\}_{n\geq 0}$ and its companion sequence $\{L_n\}_{n\geq 0}$ (see \cite{BF1} and \cite{FL1}). In \cite{Dajune} and \cite{GL}, it was shown that there are no numbers with the Lehmer property in the sequence of Cullen numbers $\{C_n\}_{n\ge 1}$ of general term $C_n=n2^n+1$, and in some appropriate generalization of the sequence of Cullen numbers, respectively.

Here, we use the same argument as in \cite{FL1} to show that there is no number with the Lehmer property in the Pell sequence $\{P_n\}_{n\ge 0}$ given by $P_0=0, P_1=1$ and $P_{n+2}=2P_{n+1}+P_n$ for all $n\geq0$.
As with other Lucas sequences, the Pell sequence has a companion $\{Q_n\}_{n\ge 0}$  given by $Q_0=2, Q_1=2$ and $Q_{n+2}=2Q_{n+1}+Q_n$ for all $n\geq 0$. There are no numbers with the Lehmer property in
the companion sequence $\{Q_n\}_{n\ge 1}$ either, but this trivially follows from the fact $Q_n$ is even for all $n\ge 1$.

Our result is the following:

\begin{theorem}
There is no Pell number with the Lehmer property.
\end{theorem}
 
 \section{Preliminary results}
Let $(\alpha,\beta)=(1+{\sqrt{2}},1-{\sqrt{2}})$ be the roots of the characteristic equation $x^2-2x-1=0$ of the Pell sequence $\{P_n\}_{n\ge 0}$. The Binet formula for $P_n$ is
\begin{equation}
\label{eq:BinetP}
P_n= \frac{\alpha^n - \beta^n}{\alpha-\beta} \quad {\text{\rm for~ all}}\quad  n\ge 0.
\end{equation}
This implies easily that the inequality
\begin{equation}
\label{eq:sizePn}
P_n\ge 2^{n/2}
\end{equation}
hold for all $n\ge 2$. Additionally, the Binet formula for $Q_n$ is
\begin{equation}
\label{eq:BinetQ}
Q_n= \alpha^n + \beta^n\quad {\text{\rm for~ all}}\quad n\ge 0.
\end{equation}
There are several relations among Pell and Pell-Lucas numbers which are well-known and can be proved using the Binet formula \eqref{eq:BinetP} for the Pell numbers and its analog \eqref{eq:BinetQ}. We only use the following well-known results.
\begin{lemma}
\label{lem:PQ}
The relation
\begin{equation}
\label{eq:relPQ}
Q_n^2 - 8P_n^2=4(-1)^n
\end{equation}
holds for all $n\ge 0$. Further, if $n$ is odd, then 
\begin{equation}
\label{eq:RelP}
P_n-1=\left\{ \begin{matrix}
P_{(n-1)/2}Q_{(n+1)/2} & {\text{if}} & n\equiv 1\pmod 4;\\
P_{(n+1)/2}Q_{(n-1)/2} & {\text{if}} & n\equiv 3\pmod 4.
\end{matrix}\right.
\end{equation}
\end{lemma}
For a prime $p$ and a nonzero integer $m$ let $\nu_p(m)$ be the exponent with which $p$ appears in the prime factorization of $m$.  The following result is well-known and easy to prove.
\begin{lemma}
\label{lem:orderof2}
The relations
\begin{itemize}
\item[(i)] $\nu_2(Q_n)=1$,
\item[(ii)] $\nu_2(P_n)=\nu_2(n)$
\end{itemize}
hold for all positive integers $n$.
\end{lemma}

\section{Proof of the Theorem}
Let us recall that if $N$ has the Lehmer property, then $N$ has to be odd and square-free.
In particular, if $P_n$ has the Lehmer property for some positive integer $n$, then Lemma \ref{lem:orderof2} (i) shows that $n$ is odd.  One checks with the computer that there is 
no number $P_n$ with the Lehmer property with $n\leq 200$. So, we can assume that $n>200$. Put $K=\omega(P_n)\geq 15$. 

From relation \eqref{eq:RelP}, we have that
$$P_n-1=P_{(n-\epsilon)/2} Q_{(n+\epsilon)/2} \quad \hbox{where}\quad \epsilon\in\{\pm 1\}.$$
By Theorem 4 in \cite{CP}, we have that $P_n<K^{2^K}.$ By \eqref{eq:sizePn}, we have that $K^{2^K}>P_n>2^{n/2}.$ Thus,
\begin{equation}
\label{eq:2}
2^K\log K>\frac{n\log 2}{2}>\frac{n}{3}.
\end{equation}
We now check that the above inequality implies that  
\begin{equation}
\label{eq:3}
2^K>\frac{n}{4\log\log n}.
\end{equation}
Indeed,  assume that  the reverse inequality \eqref{eq:3}  holds. Then,
$$
2^K\leq\frac{n}{4\log\log n},
$$
giving 
$$
K\log 2\leq \log n-\log 4 -\log\log\log n <\log n.
$$
For the above right--hand side inequality, we used the fact that $n>200>e^e$, which implies that $\log\log\log n$ is positive. Thus, 
$$
K<\frac{\log n}{\log 2}<2\log n.
$$ 
We now get that
$$
2^K\log K<\frac{n\log(2\log n)}{4\log\log n}.
$$
Using inequality \eqref{eq:2}, we have that
$$
\frac{n}{3}<2^K\log K<\frac{n\log 2}{4\log\log n} +\frac{n}{4},
$$
leading to $\log\log n<3\log 2$, which implies that $n<e^8<3000.$
However, since $K\geq 15$ and $P_n$ is odd, we have,  by Lemma \ref{lem:PQ}, that
$$
2^{15}\mid 2^{K}\mid \phi(P_n)\mid P_n-1=P_{(n-\epsilon)/2} Q_{(n+\epsilon)/2}.
$$
We observe that from Lemma \ref{lem:orderof2}, $Q_n$ is never divisible by $4$, so 
$$
2^{14}\mid P_{(n+1)/2} \quad \hbox{or}\quad 2^{14}\mid P_{(n-1)/2}.
$$
Lemma \ref{lem:orderof2} again implies that  $2^{14}=16384$ divides one of $(n+1)/2$ or $(n-1)/2$. This contradicts the fact that $n<3000$. Thus, inequality \eqref{eq:3} holds.
Let $q$ be any prime factor of $P_n$. Reducing relation
\begin{equation}
\label{eq:nodd}
Q_n^2 - 8P_n^2 = 4 (-1)^n
\end{equation}
of Lemma \ref{lem:PQ} modulo $q$, we get $Q_n^2\equiv -4\pmod q$. Since $q$ is odd, (because $n$ is odd), we get that $q\equiv 1\pmod 4$. This is true for all prime factors
$q$ of $P_n$. Hence,
$$
2^{2K}\mid \phi(P_n)\mid P_n-1=P_{(n-\epsilon)/2} Q_{(n+\epsilon)/2}.
$$
Since $Q_n$ is never divisible by $4$, we have that $2^{2K-1}\mid $ divides one $P_{(n+1)/2}$ or $P_{(n-1)/2}$. Hence, $2^{2K-1}$  divides one of $(n+1)/2$ or $(n-1)/2$. Using relation \eqref{eq:3}, we have that 
$$
\frac{n+1}{2}\geq 2^{2K-1} \geq \frac{1}{2}\(\frac{n}{4\log\log n}\)^2.
$$
This last inequality leads to 
$$
n^2<16(n+1)(\log\log n)^2,
$$
giving that $n<21$, a contradiction, which completes the proof of this theorem.

\section*{Acknowledgments} B. F. thanks OWSD and Sida (Swedish International Development Cooperation Agency) for a scholarship during her Ph.D. studies at Wits.

\end{document}